\documentstyle[a4,11pt,righttag,amssymb,verbatim]{amsart}
\newtheorem{Th}{Theorem}[section]

\newtheorem{Lem}[Th]{Lemma}
\newtheorem{Prop}[Th]{Proposition}
\newtheorem{Clm}[Th]{Claim}

\newtheorem{th}{Theorem} 
\newtheorem{prop}{Proposition}

\numberwithin{equation}{section}
\renewcommand{\theequation}{\thesection.\arabic{equation}}

\def\cB{{\cal B}} \def\cC{{\cal C}} 
  \def\Inn#1{\operatorname{Inn}(#1)}
\def\Fix{\operatorname{Fix}} \def\Cen{\operatorname{Cen}}

\def\f{\varphi} \def\s{\sigma} \def\eps{\varepsilon} 
\def\aut#1{\operatorname{Aut}(#1)}
\def\str#1{\langle #1\rangle} \def\rank#1{\operatorname{rank}#1}
\def\av#1{\overline{#1}}     \def\avst#1{\overline{\mathstrut #1}}
\def\inv{{}^{-1}} \def\id{\operatorname{id}}
\def\Z{\boldsymbol{\text Z}} 
\def\Fab{\lefteqn{\phantom{.}\overline{\phantom{I}}}{F}}
\def\fF{{\frak F}} \def\LII{{\bold L}_2} \def\strX{ X_{\text{\rm II}} }
\def\Theo{\operatorname{Th}} \def\vk{\varkappa}
\def\End#1{\operatorname{End}(#1)}   \def\gl#1{\operatorname{GL}(#1)}

\begin{document}
\title[Set theory is interpretable...]{Set theory is interpretable in the 
automorphism group
of an infinitely generated free group}
\thanks{Supported by Russian Foundation of Fundamental Research
Grant 96-01-00456}
\author{Vladimir Tolstykh}
\address{Department of mathematics, Kemerovo State University,
Krasnaja, 6,  650043, Kemerovo,  Russia}
\email{vlad\symbol{64}accord.kuzb-fin.ru}
\date{\today}
\maketitle

\section*{\it Introduction}

In 1976 S.Shelah \cite{Sh3} showed that in the endomorphism
semi-group of an infinitely generated algebra which is
free in a variety one can interpret
some set theory. It follows from his
results that, for an algebra $F_\vk$ which is free of
infinite rank $\vk$ in a variety of algebras in a language
$L,$ if $\vk > |L|,$ then the first-order theory of the endomorphism semi-group
of $F_\vk,$ $\Theo(\End{F_\vk})$ syntactically interprets
$\Theo(\vk,\LII),$ the second-order theory of the
cardinal $\vk.$  This means that for any second-order sentence $\chi$
of empty language there exists $\chi^*,$ a first-order sentence of semi-group 
language,
such that for any infinite cardinal $\vk$
$$
\chi \in \Theo(\vk,\LII) \iff \chi^* \in \Theo(\End{F_\vk}).
$$

In his paper Shelah notes that it is natural to study
the similar problem for {\it automorphism groups}
instead of endomorphism semi-groups; {\it a priori}
the expressive power of the first-order logic for
automorphism groups is less than the one for endomorphism
semi-groups. For instance, according to Shelah's results of 1973
\cite{Sh1,Sh2} on permutation groups, one {\it cannot}
interpret set theory by means of first-order logic
in the permutation group of an infinite set,
the automorphism group of an algebra in empty
language. On the other hand, one {\it can} do this in
the endomorphism semi-group of such an algebra.

In \cite{Tolsty'sPh} the author found a solution for the case of the variety
of vector spaces over a fixed division rings. If $V$
is a vector space over $D$ of an infinite dimension
$\vk,$ then theory $\Theo(\vk,\LII)$ is interpretable
in the first-order theory of $\gl V,$ the automorphism
group of $V.$ When a division ring $D$ is a countable and
definable up to isomorphism by a second-order
sentence, then theories $\Theo(\gl V))$ and $\Theo(\vk,\LII)$
are mutually syntactically interpretable.
In general case, the formulation is a bit more complicated.

The main result of this paper, Theorem \ref{MainTh}, states that
a similar result holds for the variety of all groups:

{\sc Theorem \ref{MainTh}.} {\it Let $F$ be an
infinitely generated free group of the rank $\vk.$
Then  the second-order theory of the set $\vk$ and the
elementary theory of $\aut F$ are mutually
interpretable, uniformly in $F$.}

As a corollary we have

{\sc Theorem \ref{ElEquivCrit}.} {\it
Let $F$ and $F'$ be infinitely generated free
groups of ranks $\vk$ and $\vk',$ respectively.
Then their automorphism groups are elementarily
equivalent if and only if the cardinals $\vk$ and $\vk'$ are
equivalent in the second-order logic as sets:
$$
\aut F\equiv \aut{F'} \iff \vk \equiv_{\LII} \vk'.
$$
}

Let $F$ denote an infinitely generated free group.
In Sections \ref{invs}-\ref{conjs} we
prepare `building materials' for the first-order
interpretation of infinitely generated free groups in
their automorphism groups. We prove that {\it the
subgroup $\Inn F$ of all conjugations {\rm(}inner
automorphisms of $F${\rm)} is $\varnothing$-definable in
the group $\aut F$} (Theorem \ref{DefOfConjs}). Our
key technical result states the set of all
conjugations by powers of primitive elements is
$\varnothing$-definable in $\aut F.$ (Proposition
\ref{DefOfConjsByPPE}).

We use in the proof of Theorem \ref{DefOfConjs} a
characterization of involutions in automorphism groups of
free groups -- known due to results of J.~Dyer and
G.~P.~Scott \cite{DSc}. The Theorem enables us to
prove the completeness of the automorphism groups of
arbitrary non-abelian free groups.  This generalizes
the result of J.~Dyer and E.~Formanek of 1975
\cite{DFo}, who proved  the completeness of
automorphism groups of {\it finitely} generated
non-abelian free groups.

In Section \ref{a_f_ff} we reconstruct (without
parameters) in the group $\aut F$ the three-sorted structure
$\str{\aut{F},F,S},$ where $S$ denotes the set of all
free factors of $F.$ The basic relations of the latter
structure are those of $\aut F$ and $F,$ the actions
of $\aut F$ on $F$ and $S,$ the membership relation on
$F \cup S,$ and the relation $A=B*C$ on $S.$

The main result of the Section \ref{basis} is
a recovering of a basis of $F$ in the structure
$\str{\aut{F},F,S}.$ We interpret with definable
parameters in the structure $\str{\aut F,F,S}$ the
structure $\str{\aut F,F,\cB}$ (with natural
relations), where $\cB$ is a free basis of $F.$

In the final section, using quite standard techniques,
we prove that the elementary theory of the structure
$\str{\aut F,F,\cB}$ and the second-order theory
of the set $\vk,$ where $\vk=\rank F,$ are mutually
syntactically interpretable, uniformly in $\vk.$

The author is grateful to his colleagues in Kemerovo
State University Oleg Belegradek, Valery Mishkin,
and Peter Biryukov for reading of a very draft of
this paper and helpful comments. Some results
of the paper were announced in the abstract
\cite{To}

\section{\it Involutions} \label{invs}

In what follows $F$ (unless otherwise stated) stands for an infinitely
generated free group. Let $\Fab=F/[F,F].$ Clearly,
$\Fab$ is a free abelian group of the same rank. The
natural homomorphism $w \mapsto \av w$ from the group
$F$ to $\Fab$ provides the homomorphism $\aut F \to \aut
\Fab.$ To denote this homomorphism we shall be using the
same symbol $\av{\phantom{a}}.$ We shall also say that
an automorphism $\f \in \aut F$ {\it induces} the
automorphism $\av\f \in \aut \Fab.$

In \cite{DSc} J.~Dyer and G.~P.~Scott obtained a description
of automorphisms of $F$ of prime order. For involutions
that description yields the following

\begin{Th} \mbox{\rm \cite[p. 199]{DSc} } \label{CanForm}
For every involution $\f$ in
the group $\aut F$ there is a basis $\cB$ of $F$ of the form
$$
\{u : u \in U \} \cup \{z, z' : z \in Z \} \cup
\{ x, y : x \in X, y \in Y_{x} \}
$$
on which $\f$ acts as follows
\begin{align} \label{eqCanForm}
&\quad \f u = u, \quad u \in U, \\
& \begin{cases}
\f z = z' & z \in Z, \\
\f z'= z, & \\
\end{cases}  \nonumber \\
&\begin{cases}
\f x =x \inv,   & x \in X,\\
\f y =x y x\inv, & y \in Y_x.
\end{cases} \nonumber
\end{align}
Specifically, the fixed point subgroup of $\f,$ $\Fix(\f),$ is the
subgroup $\str{u : u \in U},$ and hence is a
free factor of $F.$
\end{Th}

We shall call a basis of $F$ on which $\f$ acts
similar to (\ref{eqCanForm}) a {\it canonical} basis
for $\f.$ In view of (\ref{eqCanForm}) one can
partition every canonical basis $\cB$ for $\f$ as follows
\begin{equation} \label{eqBDivision}
\cB=U(\cB) \cup Z(\cB) \cup \{z' : z \in Z(\cB)\}\cup X(\cB) \cup
\bigcup_{x \in X(\cB)} Y_x(\cB).
\end{equation}
We shall also call any set of the form $\{x\} \cup Y_x,$
where $x \in X(\cB)$ a {\it block} of $\cB,$ and
the cardinal $|Y_x|+1$ the {\it size} of a block.
The subgroup generated by the set $Y_x$ will be
denoted by $C_x$ ($\f$ operates on this subgroup as
conjugation by $x$), and the subgroup generated by the block
$\{x\} \cup Y_x$ will be denoted by $H_x.$ Sometimes
we shall be using more `accurate' notation like
$C^{\f}_x$ or $H^{\f}_x.$ The set $U(\cB)$ will be
called the {\it fixed part} of $\cB.$

Clearly, if $\cB$ and $\cC$ are some canonical
bases for involutions $\f,\psi,$ respectively, and
the action of $\f$ on $\cB$ is isomorphic to the
action of $\psi$ on $\cC$ (that is the corresponding
parts of their canonical bases given by (\ref{eqBDivision})
are equipotent), symbolically $\f|\cB \cong \psi|\cC,$ then
$\f$ and $\psi$ are conjugate.

For the sake of simplicity
we prove the converse (in fact a stronger result)
only for involutions we essentially use: for
involutions with $Z(\cB)=\varnothing$ in all
canonical bases $\cB$'s.  We shall call these
involutions {\it soft} involutions.  It is useful that
involutions in $\aut \Fab$ induced by them have a sum of
eigen $\pm$-subgroups equal to $\Fab,$ like
involutions in general linear groups over division rings of
characteristic $\ne 2$.

We shall say that involutions $\f,\psi \in \aut F$ {\it have
the same canonical form,} if for all canonical bases
$\cB,\cC$ of $\f$ and $\psi,$ respectively, $\f|\cB
\cong \psi|\cC.$ Note that a priori we cannot even claim
that the relation we introduce is reflexive.

\begin{Prop} \label{ConjCriterion}
Let $\f \in \aut F$ be a soft involution.  An
involution $\psi \in \aut F$ is conjugate to $\f$ if and only if
$\psi$ is soft and $\f,\psi$ have the same canonical
form.
\end{Prop}

\begin{pf} Let $\Fab_2 = \Fab/2\Fab$ that is the quotient group of $\Fab$
by the subgroup of even elements. Natural homomorphisms
$F \to \Fab$ and $\Fab \to \Fab_2,$ gives us a homomorphism
$\mu : \aut F \to \aut{\Fab_2}.$ Clearly, the family of all involutions
in $\ker \mu$ coincides with the
family of all soft involutions. Therefore an involution
which is conjugate to a soft involution is soft, too.

\begin{Lem} \label{TheyMoveToInverse}
Let $\f$ be a soft involution with a canonical basis $\cB.$

{\rm (i)} Suppose  $a$ is an element in $F$ such
that $\f a=a\inv.$ Then $a=\f(w)w\inv$ or $a=\f(w) x
w\inv$ for some $x \in X=X(\cB)$ and $w \in F.$

{\rm (ii)} Suppose  $C$ is a maximal subgroup of $F$ on which
$\f$ acts as conjugation by $x \in X(\cB)${\rm:}
$$
C=\{c \in F : \f(c)=xcx\inv\}.
$$
Then $C=C^\f_x=\str{Y_x}.$
\end{Lem}

\begin{pf} (i) By \ref{CanForm} we have
\begin{equation} \label{eqDecomp}
F = \Fix(\f) * {\prod_{x \in X}}^\ast\, H_x,
\end{equation}
where each factor is $\f$-invariant. Then
$a=a_1 \ldots a_n,$ where every $a_i,$ $i=1,\ldots,n$
is an element of a free factor in expansion
(\ref{eqDecomp}), $a_i$ and $a_{i+1}$ lie in
different factors for every $i=1,\ldots,n-1$ (that
is the sequence $a_1,a_2,\ldots,a_n$ is reduced). Hence if
$\f(a)=a\inv,$ or $\f(a_1)\ldots \f(a_n)a_1 \ldots
a_n=1$ then
$$
\f(a_n)a_1=1,\quad \f(a_{n-1})a_2=1, \quad \ldots
$$
Therefore $a=\f(w_0)w_0\inv$ or $a=\f(w_0)\, v\, w_0\inv,$
where $v \in H_x$ for some $x \in X$ and $\f(v)=v\inv.$

So let us prove, using induction on length of
a word $v$ in the basis $\cB,$ that $v=\f(w_1)w_1\inv$ or
$v=\f(w_1)x w_1\inv.$ The only words $v$ of length one
in $H_x$ with $\f v=v\inv$ are $x$ and $x\inv=\f(x)x x\inv.$

An arbitrary element $v \in H_x$ can be written
in the form
\begin{equation} \label{eqZ_In_H_x}
v=x^{k_1} y_1 x^{k_2} y_2 \ldots x^{k_m} y_m,
\end{equation}
where $y_i \in C_x,$ the elements $x^{k_1}$
and $y_m$ could be equal to 1, but any other
element is non-trivial. Since $\f$ acts
on $C_x$ as conjugation by $x$ we have
\begin{equation}  \label{eqPhiOnZ}
\f(v) =x^{-k_1+1} y_1 x^{-k_2} y_2 \ldots x^{-k_m} y_m x\inv.
\end{equation}
Suppose that $\f(v)v=1.$ We then have
$$
x^{-k_1+1} y_1 x^{-k_2} y_2 \ldots x^{-k_m} y_m x\inv x^{k_1} y_1 x^{k_2} y_2 
\ldots x^{k_m} y_m=1.
$$
Let first $y_m \ne 1.$ Then $k_1=1$ and $y_m=y_1\inv.$
Hence
$$
v=xy_1x\inv (x^{k_2+1} y_2 \ldots x^{k_m}) y_1\inv=\f(y_1) t y_1\inv.
$$
It is easy to see that $\f(t)=t\inv$ and length of $t$ is less
than length of $v.$ In the case when $y_m=1$ we
have $k_m \ne 0$ and $k_1=k_m+1.$ Therefore,
$$
v=x^{k_m} (x y_1 \ldots y_{m-1}) x^{k_m} =\f(x^{-k_m}) t x^{k_m},
$$
and we again have that $\f(t)=t\inv$ and $|t|<|v|.$

(ii) Let $\f c=xcx\inv$ and $c=c_1 c_2\ldots c_n,$ where
$c_i$ are elements in free factors from (\ref{eqDecomp})
and the sequence $c_1,c_2,\ldots,c_n$ is reduced. Suppose
that $n \ge 2.$ Due to the $\f$-invariance of our
free factors, the sequence $\f(c_1),\f(c_2),\ldots,\f(c_n)$
is also reduced and must represent the same element as the
sequence $x,c_1,c_2,\ldots,c_n,x\inv.$ It is easy
to see that it is possible if both $c_1,c_n$ lie
in $H_x.$
In particular, $n \ge 3.$ It implies
that $\f(c_1)=xc_1$ and $\f(c_n)=c_n x\inv.$
It easily follows from (\ref{eqZ_In_H_x}) and
(\ref{eqPhiOnZ}) there is no
$v \in H_x$ with $\f(v)=xv;$ it of course means
that for every $v \in H_x$ the equality
$\f(v)=vx\inv$ is also impossible, since
it is equivalent to $\f(v\inv)=x v\inv.$

Thus, if $\f(c)=xcx\inv,$ then $c \in H_x.$
By applying formulae (\ref{eqZ_In_H_x}) and
(\ref{eqPhiOnZ}), one can readily conclude
that $c$ must be in $C_x.$
\end{pf}

{\sc Remarks.} (a) Note that $a=\f(w)w\inv$ cannot be
a {\it primitive} element (i.e. a member of a basis of $F$) since
$\av a$ is an {\it even} element of $\Fab.$ Indeed, it follows
from (\ref{eqDecomp}) that
$$
\avst w =\av{w(U)} + \av{w(X)} +\av{w(Y)}.
$$
where $w(U) \in \Fix(\f),$ $w(X)$ is an element in the subgroup
generated by $X,$ and $w(Y)$ is an element in the subgroup
generated by the set $Y=\bigcup_{x \in X} Y_x.$
Hence
\begin{align*}
\avst a &= \av{\f(w)w\inv}=\av{\f(w)}-\avst w \\
      &= (\av{w(U)} - \av{w(X)} +\av{w(Y)}) - (\av{w(U)} + \av{w(X)} +\av{w(Y)}) 
= \\
      &= -2 \av{w(X)}.
\end{align*}

(b) Using a similar argument we see that
if $\f(w_1) x_1 w_1\inv=\f(w_2) x_2 w_2 \inv,$ where $x_1,x_2 \in X,$
then $x_1=x_2$ (the element $\av x_1-\av x_2$ is even
if and only if $x_1=x_2$).

Suppose now that soft involution $\psi$ is a conjugate of $\f$:
$\psi=\s\inv \f \s.$ Let
$$
\cB'=U' \cup X'  \cup \bigcup_{x' \in X'} Y'_{x'}
$$
be a canonical basis for $\psi.$ Fixed point subgroups
of $\f$ and $\psi$ are clearly isomorphic. Thus,
$|U'|=|U|.$

If $\psi\, x'=x'{}\inv,$ where $x' \in X',$ then $\f(\s
x')=(\s x')\inv.$ By \ref{TheyMoveToInverse} (i) and the
above remarks there is a unique $x \in X$ such that
\begin{equation}  \label{eqs(x')=...}
\s x'=\f(w) x w\inv
\end{equation}
The mapping $x' \mapsto x$ determined in such a way is
injective, because otherwise we can find two distinct elements
in a basis of a free abelian group group whose difference is
even.

Hence, $|X'| \le |X|$ and by symmetry $|X'|=|X|.$

We claim now that
$$
\s C'_{x'} = wC_x w\inv.
$$
It will complete the proof, because in this
case $|Y'_{x'}|=|Y_x|.$

Let $y' \in C'_{x'}$ and $b=\s y'.$ Since $\psi y' =x' y' x'{}\inv,$
then $\f b= (\s x') b (\s x')\inv.$ Therefore we have
$$
\f(w\inv b w) = \f(w\inv) \f(b) \f(w)=x(w\inv b w) x\inv.
$$
Hence by \ref{TheyMoveToInverse} (ii) $w\inv b w \in C_x.$

The equation (\ref{eqs(x')=...}) can be rewritten as
follows
$$
\s\inv x = \psi(\s\inv w\inv) x' \s\inv w.
$$
Therefore
$$
\s\inv C_x \subseteq \s\inv w\inv C'_{x'} \s\inv w,
$$
or
$$
C_x \subseteq w\inv \s C'_{x'} w,
$$
and the result follows.
\end{pf}

\begin{Prop} \label{SoftsAreDef}
The set of all soft involutions is $\varnothing$-definable
in $\aut F.$
\end{Prop}

\begin{pf}
Let us now cite the above mentioned Dyer-Scott theorem in the full
form.

\begin{th}[Dyer-Scott]
Let $F$ be free, and $\alpha$ an automorphism of
$F$ of prime order $p.$ Then
$$
F=\Fix(\alpha) * ({\prod_{i \in I}}^\ast F_i) *
({\prod_{\lambda \in \Lambda}}^\ast F_\lambda),
$$
where each factor is $\alpha$-invariant. Moreover,

{\rm (i)} For each $i\in I,$ $F_i$ has a basis
$x_{i,1}\ldots,x_{i,p}$ such that
$$
\alpha(x_{i,r})=x_{i,r+1 (\text{\rm mod}\, p)}.
$$

{\rm (ii)} For each $\lambda \in \Lambda,$ $F_\lambda$
has a basis
$$
x_{\lambda,1},\ldots,x_{\lambda,p-1},\quad \{y_j : j \in J_\lambda\}
$$
such that
\begin{alignat*}3
&\alpha(x_{\lambda,r})    &&=x_{\lambda,r+1},       &\quad &r=1,\ldots,p-2,\\
&\alpha(x_{\lambda,p-1}) &&= (x_{\lambda,1},\ldots,x_{\lambda,p-1})\inv,\\
&\alpha(y_j)             &&=x_{\lambda,1}\inv y_j x_{\lambda,1}, &\quad &j \in 
J_\lambda.
\end{alignat*}
\end{th}

Let $\s \in \aut F$ be an element of prime order $p > 2.$
Consider the natural homomorphism from $\aut F$ to $\aut {\Fab_2}.$
As an easy corollary of the Dyer-Scott Theorem, we have that
the image of $\s$ under this homomorphism is non-trivial.

Two involutions $\f,\f_0$ in the automorphism
group of two-generator free group $G=\str{a_1,a_2}$ such that
$$
\begin{cases}
\f(a_1)=a_2\inv,\\
\f(a_2)=a_1\inv,
\end{cases}
\quad
\begin{cases}
\f_0(a_1)=a_1a_2,\\
\f_0(a_2)=a_2\inv,
\end{cases}
$$
are non-soft ($\{a_1,a_2a_1\}$ is a canonical
basis for latter, and $\{a_1,a_2\inv\}$ for former).
Their product $\s=\f_0\f$ is of the order three in $\aut{G}$:
\begin{align*}
\s(a_1) &=a_2,\\
\s(a_2) &=(a_1 a_2)\inv.
\end{align*}
Generalizing this example, we can easily observe
that for every non-soft involution $\f \in \aut F$
there is a conjugate $\f_0 \in \aut F$ of $\f$
such that $\f_0\f$ has order three.

The latter (first-order) property is evidently false
for every soft involution $\f \in \aut F,$ because for
every conjugate $\f_0$ of $\f,$ the product $\f_0\f$
cannot have order three (all elements of order three
in $\aut F$ have non-trivial images under the natural
homomorphism from $\aut F$ in $\aut {\Fab_2},$ but
$\f_0\f$ has not).
\end{pf}


\section{\it Anti-commutative conjugacy classes} \label{acc}

The key roles in the group-theoretic characterization
of conjugations in $\aut F$ will be played by two conjugacy classes
of involutions. First class consists of involutions
with a canonical form
\begin{alignat*}2
\f x &=x \inv, & \quad & \\
\f y &= x y x\inv, & y &\in Y
\end{alignat*}
that is $\cB=\{x\} \cup Y$ is a basis of $F,$ canonical for $\f,$
$U(\cB)=Z(\cB)=\varnothing,$
$X(\cB)$ is a singleton set $\{x\},$ $Y_x=Y$ (any canonical
form reproduced below is interpreted in a similar way).
We shall call these involutions {\it quasi-conjugations.}

An arbitrary element $\f \in \aut F$ in the second
class has the following canonical form:
$$
\f x =x\inv,\quad x \in X,
$$
that is there is a basis of $F$ such that $\f$ inverts
all its elements. We shall use for these involutions
the term {\it symmetries.}

Every symmetry induces in $\aut \Fab$ the automorphism
$-\id_{\Fab},$ and hence the product of two symmetries
induces $\id_{\Fab}.$ Therefore by \ref{CanForm} two
symmetries commute if and only if they are coincide.
Thus, the conjugacy class of all symmetries is, say,
{\it anti-commutative}, since its elements are
pairwise non-commuting. In the next section we shall
prove that the class of all quasi-conjugation is also
anti-commutative.

In order to characterize conjugations we shall use
anti-commutative conjugacy classes of involutions, but we need not an
exact determination of {\it all} such conjugacy classes: it
suffices to know that they lie in some
`easy-to-define' family. In the following proposition
we formulate and prove a necessary condition of being
an anti-commutative conjugacy class, but do not prove its
sufficiency.

\begin{Prop} \label{ACCs}
Let $\f$ be an involution in an anti-commutative
conjugacy class. Then either $\f$
has a canonical form such that
\begin{alignat}2 \label{eqBeads}
\f u &=u,  &\quad u &\in U,\\
\f x &=x\inv,  & x & \in X,  \nonumber \\
\f y &=x y x\inv,  & y & \in Y_x, \nonumber
\end{alignat}
where $|X| \ge 2,$ all the sets $Y_x,$ $x \in X$ have
the same finite power $n$ and $|U| < n+1,$ or $\f$ has
a canonical form such that
\begin{alignat}2 \label{eqSnakes}
\f u &=u,  &\quad u &\in U,\\
\f x &=x\inv,      &\quad &  \nonumber \\
\f y &=x y x\inv,  & y & \in Y, \nonumber
\end{alignat}
where the cardinal $|U|$ is finite and less than $|Y|+1.$

Furthermore, all involutions of the form {\rm (\ref{eqBeads})}
are squares in $\aut F$ and involutions of the form
{\rm (\ref{eqSnakes})} are not.
\end{Prop}

In other words, in the terminology introduced in
Section \ref{invs}, a canonical basis $\cB$ for an
involution in an anti-commutative conjugacy class
either contains exactly one block and the power of the
fixed part of $\cB$ is less than the size of this
block, or all blocks of $\cB$ have the same finite
size and the power of the fixed part is less than the
size of any block. Clearly, symmetries have the form
\eqref{eqBeads} (all blocks of their canonical bases
have the size one), and quasi-conjugations have the
form \eqref{eqSnakes} (the size of the unique block
is equal to $\rank F$).

\begin{pf}
Show first that every anti-commutative class of
involutions consists only of soft involutions. Indeed, let
involution $\f$ be an involution whose canonical basis
$\cB$ has non-empty `permutational' part $Z(\cB).$
Suppose $\f$ takes $z \in Z(\cB)$ to $z'.$
Consider an involution $\psi$ which acts on $\cB
\setminus \{z,z'\}$ exactly as $\f$ does, but taking
$z$ to $z'{}\inv.$ Clearly, $\psi$ is conjugate to
$\f$ and commutes with $\f.$

The following example demonstrates why the size of the fixed part of
a canonical basis must necessarily be less than the size of each block:
$$
\begin{cases}
\f x=x\inv,\\
\f y=x y x \inv,\\
\f x_1=x_1,\\
\f y_1=y_1,\\
\f u=u
\end{cases}
\quad
\begin{cases}
\psi x =x,\\
\psi y =y,\\
\psi x_1 =x_1\inv,\\
\psi y_1 =x_1 y_1 x_1 \inv,\\
\psi u =u.
\end{cases}
$$

Let us make a technical remark. Involutions
\begin{equation}  \label{eqTwoMinuses}
\begin{cases}
\f x=x \inv,  \\
\f y=x y x\inv,\, y \in Y
\end{cases}
\mbox{   and   }
\begin{cases}
\psi x=x \inv, \\
\psi y=x\inv y x\inv,\, y \in Y
\end{cases}
\end{equation}
are conjugate: the second one acts `canonically' on
the set $\{x, x\inv y\ : y \in Y\}$: $\f(x\inv y)=
x(x\inv y)x\inv.$

Let now $\f$ be a soft involution with a canonical
basis $\cB$ such that for some distinct $x_1, x_2 \in
X(\cB)$ either $|Y_{x_1}| < |Y_{x_2}|$ or
$|Y_{x_1}|=|Y_{x_2}|$ and the cardinal $|Y_{x_1}|$ is
infinite (thus, there must be no neither a pair of
blocks of different size nor a pair of infinite
blocks). The fact that $\f$ commutes with at least two
its conjugates is a consequence of
(\ref{eqTwoMinuses}) and the following

\begin{Clm}
Let $G$ be a free group with a basis
$$
\{x,a\} \cup B \cup \{c\} \cup D \cup E,
$$
where $|B|=|D|.$ Then

{\rm (i)} Involutions
$$
\begin{cases}
\f x=x\inv,\\
\f a=x\inv a x\inv,\\
\f b=x\inv b x\inv,\\
\f c=c\inv,\\
\f d=c\inv d c\inv,\\
\f e=x\inv e x\inv
\end{cases}
\qquad \text{ and }\qquad
\begin{cases}
\psi x=x\inv,\\
\psi a=a\inv,\\
\psi b=a\inv b a\inv, \, b \in B,\\
\psi c=x\inv c x\inv,\\
\psi d=x\inv d x\inv, \, d \in D,\\
\psi e=x\inv e x\inv, \, e \in E
\end{cases}
$$
are conjugate and commute.

{\rm (ii)} Furthermore,
\begin{align*}
&\rank C^{\f}_x = |\{a \} \cup B \cup E| =|B|+|E|+1,\\
&\rank C^{\f}_c=|D|.
\end{align*}
Thus, $\rank C^{\f}_c \le \rank C^{\f}_x,$ where
equality holds only if both ranks are infinite.
\end{Clm}
\begin{pf}
Easy.
\end{pf}

To complete the proof of the first statement in the
Proposition, we should cut off involutions with a
canonical form such that
\begin{alignat*}2
\f u &=u, & u &\in U,\\
\f x &=x\inv, &&\\
\f y &=xyx\inv, &\quad y &\in Y,
\end{alignat*}
where $|U|$ is infinite. To do this let us
consider involutions $\f,\psi$ which act
on a basis $U \cup \{x\} \cup Y$ of $F$
as follows
$$
\begin{cases}
\f u_0=u_0,\\
\f u=u,\\
\f x=x\inv,\\
\f y=xyx\inv,
\end{cases}
\quad
\begin{cases}
\psi u=u_0\inv,\\
\psi u=u,\,        \qquad     &u \in U \setminus \{u_0\}, \\
\psi x=u_0 x u_0\inv,\\
\psi y=u_0 y u_0\inv,  \quad  &y \in Y,
\end{cases}
$$
where $u_0$ is a fixed element in $U.$ It is readily seen that
$\f$ and $\psi$ are conjugate
and commute.

One can readily check that
$$
\begin{cases}
\s x =a\inv,\\
\s y =a b a\inv,\\
\s a =x,\\
\s b =y
\end{cases}
\Longrightarrow
\begin{cases}
\s^2 x =x\inv,\\
\s^2 y =x y x\inv,\\
\s^2 a =a\inv,\\
\s^2 b   =a b a\inv.
\end{cases}
$$
It demonstrates that every automorphism of $F$ of the form (\ref{eqBeads})
is a square in $\aut F.$

On the other hand, every automorphism of the form (\ref{eqSnakes})
cannot be a square. Indeed, every  $\f\in \aut F$ of the form
(\ref{eqSnakes}) induces in $\aut \Fab$ an {\it extremal} involution,
that is an involution such that one of its $\pm$-eigen subgroups,
$\str{\av x}$ in our case, has rank one. Suppose that
$\f=\s^2$ in $\aut F.$ We then have
$$
\av{\f}(\av{\s}\,\av x)= \av{\s}(\av{\f}\,\av x)=-\av{\s}\,\av x.
$$
Therefore $\av{\s}\,\av x=m\,\av x$ for $m \in \Z.$ Since
$-1$ is not a square in $\Z,$ the equation $\av\f=\av\s^2$
is impossible.
\end{pf}

\section{\it Quasi-conjugations} \label{q-conjs}

In this section we obtain a first-order
characterization of quasi-conjugations in $\aut F.$
First we prove that the conjugacy class of all
quasi-conjugations is anti-commutative. Then we
distinguish this class from other anti-commutative
conjugacy classes of involutions:  it is trivial in
the case when $\rank F=2$ and more technical in the
case when $\rank F > 2.$

\begin{Prop} \label{QsConjs:ACC}
The class of all quasi-conjugation is an anti-commutative
conjugacy class.
\end{Prop}

\begin{pf}
Let $\f$ be a quasi-conjugation with a canonical
basis $\cB=\{x\} \cup Y$:
\begin{alignat*}2
\f x &=x \inv, & \quad & \\
\f y &= x y x\inv, & y &\in Y
\end{alignat*}
Every $\s \in \aut F$ which commutes
with $\f$ takes $x$ to a primitive element of the form
$\f(w)xw\inv,$ $w \in F$ (Lemma
\ref{TheyMoveToInverse}).  It turns out that the only
primitive elements of the form $\f(w) x w\inv$ are $x$
and $x\inv.$ This fact is a consequence of the
following result, which we shall use once more later.

\begin{Lem} \label{X_And_X_inv}
Let $\alpha$ be an involution with a canonical form
such that
\begin{alignat*}2
\alpha\, x &=x\inv, &\quad&\\
\alpha\, y &=xyx\inv, &y &\in Y=Y_x,\\
\alpha\, u &=u,       &u &\in U,
\end{alignat*}
and $a\in F$ a primitive element which $\alpha$
sends to its inverse. Then $a=vx^{\pm1} v\inv,$
where $v \in \str U =\Fix(\alpha).$
\end{Lem}

\begin{pf}
As we observed earlier $a=\alpha(w)xw\inv.$
It is easy to see that $a$ lies in the normal
closure of $x.$ One can use induction on
length of a reduced word $w$ in the basis
$\{x\} \cup Y \cup U.$ We have
\begin{align*}
&\alpha(xw_1)xw_1\inv x\inv=x\inv (\alpha(w_1)xw_1\inv) x\inv,\\
&\alpha(uw_1)xw_1\inv u\inv=u(\alpha(w_1)xw_1\inv)u\inv,\\
&\alpha(yw_1)xw_1\inv y\inv= x \cdot y(x\inv \alpha(w_1)xw_1\inv)y\inv.
\end{align*}

\begin{prop} \label{Magnus'sProp}
\mbox{\rm (\cite{Ma}, \cite[II.5.15]{LSch})}
Let $F$ be a free group and the normal closure of $q
\in F$ consists of a primitive element $p.$ Then $q$
is conjugate to $p$ or $p\inv.$
\end{prop}

Hence $a=bx^\eps b\inv,$ where $\eps=\pm 1.$ Since
$\alpha(a)=a\inv,$ we have
$$
\alpha(b) x^{-\eps} \alpha(b\inv)=b x^{-\eps} b\inv.
$$
It follows that $x^{-\eps}$ and $b\inv \alpha(b)$
commute. Therefore both these elements lie
in a cyclic subgroup of $F$ (\cite[I.2.17]{LSch}).
It must be the subgroup $\str x,$ because $x$ is primitive.
Hence $\alpha(b)=b x^k.$ If $k$ is even, say
$k=2m,$ then $\alpha(b x^m)=bx^m$ and $v=bx^m  \in \str U.$
In the case when $k$ is odd, we have  there is
$z \in F$ such that $\alpha(z)=zx.$ One can easily check
that this is impossible (e.g. using `abelian' arguments
as in Section \ref{invs}).
\end{pf}

\begin{Prop} \label{CenOfASnake}
The centralizer of a quasi-conjugation $\f$ with
a canonical basis $\{x\} \cup Y$ consists of automorphisms $\s$ of $F$
of the form
\begin{alignat}2  \label{eqCenOfASnake1}
\s x &=x,  &\quad & \\
\s y &= \theta(y), & y &\in Y,  \nonumber
\end{alignat}
and of the form
\begin{alignat}2 \label{eqCenOfASnake2}
\s x &=x\inv,  &\quad & \\
\s y &=x \theta(y) x\inv, & y &\in Y, \nonumber
\end{alignat}
where $\theta \in \aut{C_x}.$
\end{Prop}

\begin{pf}
The Proposition easily follows from Lemma \ref{X_And_X_inv}
and one more

\begin{Lem} \label{A*B}
\mbox{\rm \cite[p. 101]{Cohen}}
Let $F=G*H$ be a free factorization of a free group
$F.$ Suppose that $\alpha \in \aut F$ and $\alpha|_{G}$ is an
endomorphism of $G.$ Then $\alpha|_{G} \in \aut G$ if
either $\rank G$ or $\rank H$ is finite.
\end{Lem}

By \ref{X_And_X_inv} if $\s \in \Cen(\f),$ then either
$\s x=x$ or $\s x=x\inv.$ Assume that $\s$ fixes $x.$
Then, for each $y \in C_x$
$$
\f\s y=\s x\, \s y\, \s x\inv= x\, \s y\, x\inv.
$$
Therefore $\s y \in C_x.$

If $\s x=x\inv,$ then $\f\s x=x.$ It follows that $\f\s$ has
the form (\ref{eqCenOfASnake1}), and hence
$\s$ must have the form (\ref{eqCenOfASnake2}).
\end{pf}

We can prove now that $\f$ is the unique
quasi-conjugation in its centralizer, or, in other
words, the conjugacy class of $\f$ is
anti-commutative.  Indeed, there are no
quasi-conjugations in the family of automorphisms of
the form (\ref{eqCenOfASnake1}), because every
quasi-conjugation has trivial fixed-point subgroup.
Let $\s$ be a quasi-conjugation of the form
(\ref{eqCenOfASnake2}). Hence $\theta^2=\id,$ and
$\theta$  is either the identity automorphism of $C_x$ or a
soft involution. Assume that $\theta$ is an
involution. By Theorem \ref{CanForm} there is a
canonical basis $\cC$ of a free group $C_x$ for
$\theta$ such that
\begin{alignat*}3
&\theta & u &=u,              & u &\in U(\cC),\\
&\theta & c &=c\inv,          & c &\in X(\cC),\\
&\theta & d &=c d c\inv,\quad & d &\in Y_c(\cC).
\end{alignat*}
Therefore
\begin{alignat*}3
&\s  x    &&=x\inv,                                  \\
&\s  u    &&=x u x\inv,              & u &\in U(\cC),\\
&\s  (xc) &&=(xc)\inv,               & c &\in X(\cC),\\
&\s  d    &&=(xc) d (xc)\inv,\quad   & d &\in Y_c(\cC).
\end{alignat*}
We obtain a canonical basis for $\s.$ This basis
contains at least two blocks, and hence by Proposition
\ref{ConjCriterion} $\s$ cannot be a quasi-conjugation.
So $\theta$ must be the identity, or, equivalently,
$\s=\f$ as desired. This completes the proof. \end{pf}

The next is

\begin{Th}
The class of all quasi-conjugations is the unique
anti-commutative conjugacy class $K$ of involutions
such that its elements are not squares and for every
anti-commutative conjugacy class $K'$ of involutions,
whose elements are squares, all involutions in $KK'$
are conjugate.
\end{Th}

\begin{pf} Let $\f$ be a quasi-conjugation, $K'$
an anti-commutative conjugacy class of involutions, and
$\f \not\in K'.$ It suffices to prove that if
$\psi,\psi' \in K'$ both commute with $\f,$ then
$\f\psi$ and $\f\psi'$ are conjugate.

Suppose first that fixed point
subgroups of $\psi$ and $\psi'$ are trivial. Then both
$\psi$ and $\psi'$ have the form
\eqref{eqCenOfASnake2}:
\begin{alignat*}2
\psi x &=x\inv,                  & \psi' x &=x\inv \\
\psi y &=x \theta(y) x\inv,\quad & \psi' y &=x\theta'(y) x\inv, \quad y \in Y,
\end{alignat*}
where $\theta$ and $\theta'$ are in $\aut{C_x}.$ It is easy to see that
$\theta$ and $\theta'$ are soft involutions. Let $\cC$ be a canonical
basis for $\theta.$ As we observed above, the set
\begin{equation}
(\{x\} \cup U(\cC)) \cup \bigcup_{c \in X(\cC)} (\{xc\} \cup Y_c(\cC))
\end{equation}
is a canonical basis for $\psi,$ and (\theequation) is a partition
of this basis into blocks (there are at least two blocks,
because $\psi$ must be a square). Since $\psi$ lies in
an anti-commutative conjugacy class, all these blocks have the same
size, and hence for all $c \in X(\cC)$
$$
|U(\cC)|=|Y_c(\cC)|.
$$
By applying a similar argument to $\psi',$ we see
that $\theta$ and $\theta'$ are conjugate in $\aut {C_x}.$
Hence $\f \psi$ and $\f \psi'$ are conjugate in $\aut F$:
both these automorphisms fix $x$ and their restrictions
on $C_x$ ($\theta$ and $\theta'$) are conjugate.

Assume that $\psi$ and $\psi'$ have non-trivial fixed
point subgroups. Both $\psi$ and $\psi'$ preserve the
subgroup $C_x$ and fix the element $x.$ It means
(Prop. \ref{CanForm}, Prop. \ref{ConjCriterion}) that restrictions of
$\psi$ and $\psi'$ on $C_x,$ say, $\theta$ and
$\theta',$ respectively are conjugate in $\aut{C_x}$:
$\theta'=\pi\inv \theta \pi,$ where $\pi \in \aut{C_x}.$
The map $\pi$ can be extended to an element
$\s \in \aut F$ such that $\s x=x.$ By \ref{CenOfASnake}
$\s$ commutes with $\f,$ and hence
$$
\s\inv(\f \psi)\s=\f \psi'.
$$

Let us prove the converse. It suffices to prove that
for each involution $\f$ of the form \eqref{eqSnakes},
which is not a quasi-conjugation, there are two {\it
symmetries} $\psi$ and $\psi'$ such that $\f\psi$ and
$\f\psi'$ are non-conjugate involutions.

(a) A natural way to define an action of a symmetry commuting
with $\f$ on a block of a canonical basis for $\f$
is given by Proposition \ref{CenOfASnake}:
\begin{alignat*}2
\f x &=x\inv,           & \psi x &=x\inv,\\
\f y &=x y x\inv \quad  & \psi y &=x y\inv x\inv,\quad (\psi(xy)=(xy)\inv) \quad 
y \in Y_x.
\end{alignat*}
Clearly, the product of $\f$ and $\psi$ fixes $x$ and invertes
each element in $Y_x.$

(b) A way to define an action of a symmetry on the fixed
part of a canonical basis is obvious: a symmetry should invert
each element.

Let now $\cB=U \cup \{x \} \cup Y$ be a canonical basis for $\f$:
\begin{alignat*}2
\f u &=u, &\quad u &\in U \quad (U \ne \varnothing),\\
\f x &=x\inv, && \\
\f y &=xyx\inv, & y &\in Y
\end{alignat*}
One can easily find a symmetry $\psi$ such that the product
$\f\psi$ will be non-conjugate to each product of $\f$ with a symmetry
obtained in a natural way, using (a) and (b):
\begin{alignat*}4
&\psi && u_0  &&=u_0\inv, &&\\
&\psi && u    &&=u\inv, &  u &\in U \setminus \{u_0\},\\
&\psi && x    &&=u_0 x\inv u_0\inv, &&\\
&\psi && y    &&=u_0 xy\inv x\inv u_0\inv, &\quad y &\in Y,
\end{alignat*}
where $u_0 \in U.$ The reason is that any canonical
basis for $\f\psi$ contains a block of the size two
($\{u_0\} \cup \{x\}$ in our example).

The proof of Theorem \theTh\ is now complete.
\end{pf}

All the hypotheses in the latter Theorem are surely first-order,
and hence we have

\begin{Th} \label{DefOfQsConjs}
The set of all quasi-conjugations is first-order definable
in $\aut F.$
\end{Th}

{\sc Remark.} In fact, the set of all quasi-conjugations
is definable in the automorphism group of an {\it arbitrary}
non-abelian free group (whether finitely or infinitely
generated). We does not prove this result here, but discuss
key technical points. It is easy to see that the method
of the proof of Proposition \ref{QsConjs:ACC} works
for arbitrary free groups of rank at least two. Then
there are no further problems in a proof of definability
of quasi-conjugations in automorphisms groups of two-generator
free groups:

\begin{prop}
Suppose that $\rank F=2.$ Then
the class of all quasi-con\-ju\-ga\-tions is the
unique anti-com\-mu\-ta\-tive conjugacy class $K$ of
involutions whose elements are not squares.
\end{prop}

However, the property of being of a square
does not always distinguish involutions of the form
\eqref{eqBeads} from involutions of the form
\eqref{eqSnakes}, and the generalization of Theorem
\ref{DefOfQsConjs} is formulated as follows:

\begin{th}
Let $\rank F > 2.$ Then the class of all quasi-conjugations is the
unique anti-commutative conjugacy class $K$ of
involutions such that for every anti-com\-mu\-ta\-tive
conjugacy class $K'$ of involutions, all involutions in $KK'$ are
conjugate.
\end{th}

\section{\it Conjugations} \label{conjs}

\begin{Th} \label{DefOfConjs}
The subgroup $\Inn F$ of all conjugations is
$\varnothing$-definable subgroup of $\aut F.$
\end{Th}

\begin{pf} We start with a quasi-conjugation $\f.$
Suppose $\cB=\{x\} \cup Y$ is a canonical basis for
$\f.$ Let $\Pi$ denote the set of all automorphisms of $F$
of the form $\pi=\s\s',$ where $\s$ and $\s'$ are in $\Cen(\f)$
and conjugate. By \ref{CenOfASnake} conjugate
$\s,\s'$ in the centralizer of $\f$ either both have
the form \eqref{eqCenOfASnake1} (when their fixed point
subgroups are non-trivial) or have the form \eqref{eqCenOfASnake2}.
Therefore every $\pi \in \Pi$ has the form
\begin{align*}
\pi x &=x\\
\pi y &=\theta(y),\quad y \in Y,
\end{align*}
where $\theta \in \aut{C_x},$ that is $\pi$ fixes $x$ and
preserves the subgroup $C_x.$

In analogy to the term in the linear group theory
we call any involution $\psi \in \aut F$ with a canonical form such that
\begin{align*}
\psi x &=x\inv,\\
\psi y &=y, \quad y \in Y
\end{align*}
{\it extremal} involution (compare with Section
\ref{acc}).  Actually we need not extremal involutions
in the proof of the Theorem, but we shall use them
later.

\begin{Lem} \label{AllWeNeedInCen(Pi)}
\mbox{\rm (i)} All conjugations by powers of $x$ are
in the centralizer of the family $\Pi.$ Every
member of $\Cen(\Pi)$ is either an involution
or conjugation by a power of $x.$

\mbox{\rm (ii)} Every involution in $\Cen(\Pi)$ is either
a quasi-conjugation or an extremal involution. Therefore
the set of extremal involutions is $\varnothing$-definable
in $\aut F.$
\end{Lem}

\begin{pf} Let
$$
\cC=\{a,b\} \cup C
$$
be a basis of $C_x,$ and $\tau \in \aut F$ an element
of $\Cen(\Pi).$

First we construct $\pi \in \Pi$ such that the
fixed point subgroup of $\pi$ is the subgroup
$\str{x,a}.$ Since $\tau$ must commute with $\pi$ we shall have that
$$
\tau a=w_a(x,a),
$$
where $w_a$ is a reduced word in letters $x$ and
$a.$

To construct $\pi,$ let us consider conjugate $\s,\s'$ in $\Cen(\f)$
which act on $\cC$ as follows
\begin{alignat*}2
\s x &=x,                   & \s' x &=x, \\
\s a &=a\inv,               & \s' a &=a\inv,\\
\s b &=b\inv,               & \s' b &=a\inv b\inv a,\\
\s c &=c\inv, \qquad        & \s' c &=a\inv c\inv a,\quad c \in C.
\end{alignat*}
The restriction of $\pi=\s\s'$ on $C_x$ is conjugation
by $a.$ Then it is easy to show that the fixed point
subgroup of $\pi$ is $\str{x,a}.$ By Lemma \ref{A*B}
$\tau \str{x,a}=\str{x,a}.$

A similar argument can be applied to an arbitrary
primitive element in $C_x.$ Hence for every primitive
$d \in C_x$
$$
\tau d=w_d(x,d),
$$
and $\tau$ preserves the subgroup $\str{x,d}.$ We then have
$$
\tau \str x=\tau(\str{x,a} \cap \str{x,b})=\str{x,a} \cap \str{x,b}=\str x.
$$
Therefore $\tau x=x^{\pm 1}.$ In particular, the word $w_a(x,a)$
must have explicit occurrences of $a.$

We claim now that the words $w_d,$ where $d=a,b,ab$ have the same structure,
that is any word $w_d(x,d)$ can be obtained from the word
$w_a(x,a)$ by replacing occurrences of $a$ by $d$:
$$
[w_a(x,a)]^a_d=w_d(x,d).
$$
To prove this, it suffices to find in $\Pi$ an automorphism
which takes $a$ to $b$ ($a$ to $ab$).

Let $\s_1$ and $\s'_1$ be involutions in $\Cen(\f)$
such that $\s_1$ and $\s'_1$ both fix
the set $\{x\} \cup C$ pointwise and
\begin{alignat*}2
\s_1 a &=b\inv,          & \s'_1 a &=ab, \\
\s_1 b &=a\inv, \qquad   & \s'_1 b &=b\inv,\\
\end{alignat*}
Clearly, $\s_1$ and $\s'_1$ are conjugate and
$\pi_1=\s_1\s'_1$ sends $a$ to $b.$ Since $\tau$
and $\pi_1$ commute, we have
\begin{align*}
\tau a=w_a(x,a) &\Rightarrow \tau (\pi_1 a) =w_a(\pi_1 x,\pi_1 a) \\
		&\Rightarrow w_b(x,b)=w_a(x,b).
\end{align*}
Thus, there is a reduced word $w$ in letters $x$ and, say, $t$ such that
$$
[w(x,t)]^t_d=w_d(x,d),
$$
where $d=a,b,ab.$ We then have
$$
\tau (ab)=w(x,ab)=\tau(a) \tau (b)=w(x,a)w(x,b),
$$
and hence
\begin{equation}
w(x,ab)=w(x,a) w(x,b).
\end{equation}

Now we show that the word $w(x,t)$ has the form
$x^k t x^{-k},$ where $k \in \Z.$ Assume that $w(x,t)$ has
the (possibly non-reduced) form
such that
$$
x^{k_1} t^{l_1} x^{k_2} t^{l_2} \ldots x^{k_m} t^{l_m},
$$
where $k_1$ or $l_m$ could be equal to zero, whereas any other
exponent is non-trivial. Then by (\theequation)
\begin{equation}
x^{k_1} (ab)^{l_1} x^{k_2} (ab)^{l_2} \ldots x^{k_m} (ab)^{l_m} =
x^{k_1} a^{l_1} x^{k_2} a^{l_2} \ldots x^{k_m} a^{l_m}
            x^{k_1} b^{l_1} x^{k_2} b^{l_2} \ldots x^{k_m} b^{l_m}.
\end{equation}
The latter equality is evidently impossible when $m \ge 2$ and $l_m \ne 0.$
Hence $l_m=0$ and $k_m=-k_1.$ Even after this reduction (\theequation)
fails, if $m \ge 3.$ Therefore
$$
x^{k_1} (ab)^{l_1} x^{-k_1}=x^{k_1} a^{l_1} b^{l_1} x^{-k_1},
$$
and we have
$$
(ab)^{l_1}=a^{l_1} b^{l_1}.
$$
Since $a$ and $b$ are independent, $l_1=1.$

Summing up, we see that $\tau$ acts on $\cB$ as follows
\begin{align}
\tau x &=x^\varepsilon,  \\
\tau a &=x^k a x^{-k}, \nonumber  \\
\tau b &=x^k b x^{-k},   \nonumber\\
\tau c &=x^k c x^{-k}, \quad c \in C. \nonumber
\end{align}
where $\varepsilon=\pm 1.$ In the case when
$\varepsilon=-1,$ $\tau$ is an involution, otherwise
$\tau$ is conjugation by $x^k.$ Conversely, every automorphism of $F$
of the form (\theequation) is in $\Cen(\Pi).$

Suppose $\tau x=x\inv.$ Then $\tau a=x^{2m} a x^{-2m},$ where $m \in \Z$
is equivalent to
$$
\tau (x^m a x^{-m}) =x^m a x^{-m}.
$$
It demonstrates that $\tau$ of the form (\theequation)
is an extremal involution, if $\varepsilon=-1$ and $k$
is even. Clearly, if $k$ is odd, then $\tau$ is a
quasi-conjugation.
\end{pf}

Thus, we have that

\begin{Prop} \label{DefOfConjsByPPE}
The set of all conjugations by powers of primitive elements
is $\varnothing$-definable in $\aut F.$
\end{Prop}

It is easy to see now that the subgroup of all conjugations
is definable in $\aut F.$ Indeed, let $\cB$
is a basis of $F.$ Suppose that $a=w(b_1,\ldots,b_n),$
where $a$ is a non-trivial element in $F$ and $b_1,\ldots,b_n \in \cB.$
Since $\cB$ is {\it infinite}, there is $b \in \cB \setminus 
\{b_1,\ldots,b_n\}.$
Then
$$
b \text{  and  } b\inv w(b_1,\ldots,b_n)
$$
are both primitive elements. Therefore $a$ is the product of two
primitive elements. It implies that an arbitrary conjugation
in $\aut F$ is the product of two conjugations by powers of primitive
elements.
\end{pf}

{\sc Remark.} Modulo definability of
quasi-conjugations, Proposition \ref{DefOfConjsByPPE}
remains true for arbitrary non-abelian free groups. A
proof for free groups with more than two generators is
the same as the proof just completed.  We give a
sketch of proof for two-generator free groups.  Let
$\{x\} \cup \{y\}$ be a canonical basis for a
quasi-conjugation $\f.$ The unique involution in the
centralizer of $\f$ commuting with at least two its
conjugates is an involution $\psi$ which fixes $x$ and
inverts $y$ (Proposition \ref{CenOfASnake}).  The
centralizer of $\psi$ consists of conjugations by the
powers of $x,$ involutions, but some elements of
infinite order.  Fortunately, these elements can be
distinguished from conjugations: they induces in
$\aut{\Fab}$ automorphisms whose determinates
equal to $-1,$ and hence one can not represent them as
the product of conjugate elements from $\aut F.$ On
the other hand, one can represent every conjugation by
a power of $x$ as the product of two conjugate
involutions (like in the proof of Lemma \ref{AllWeNeedInCen(Pi)}).

It is known that the automorphism group of a centreless group
$G$ is complete if and only if the subgroup $\Inn G$ is a characteristic
subgroup of $\aut G$ (\cite{Spe}).
Being definable in $\aut F,$ the subgroup $\Inn F$ is
its characteristic subgroup. Therefore

\begin{Th} \label{Aut(F)IsComplete}
Let $F$ be an infinitely generated free group. Then
the group $\aut F$ is complete.
\end{Th}

As we noted in Introduction, the latter theorem
generalizes the result of J.~Dyer and E.~Formanek,
stating the completeness of automorphism groups of
finitely generated non-abelian groups \cite{DFo}.
One more purely algebraic result which is a consequence
of Theorem \ref{DefOfConjs} is the following

\begin{Th}  The automorphism groups of infinitely
generated free groups $F$ and $F'$ are isomorphic if and
only if $F \cong F'.$
\end{Th}

\begin{pf}
Any isomorphism from $\aut F$ to $\aut{F'}$ preserves
conjugations, and hence induces an isomorphism
between $F$ and $F'.$
\end{pf}

\section{\it The group and its free factors}   \label{a_f_ff}

Let us consider the structure $\fF$ such that the domain of
$\fF$ consists elements of three sorts:
\begin{itemize}
\item the elements of the group $F,$

\item the automorphisms of $F,$

\item the free factors of $F.$
\end{itemize}
The basic relations of $\fF$ are
\begin{itemize}
\item those of standard relations on sorts,

\item the actions of $\aut F$ on other sorts,

\item the membership relation between elements
of $F$ and the set of free factors;

\item a ternary relation, say $R,$ on the
set free factors such that
$$
R(A,B,C) \leftrightarrow A=B*C.
$$
\end{itemize}

\begin{Th}
The structure $\fF$ is interpretable without
parameters in $\aut F$ by means of first order logic.
\end{Th}

\begin{pf}
By Theorem \ref{DefOfConjs} the subgroup of  conjugations
is $\varnothing$-definable in $\aut F.$ Let $\tau_a$ denote
conjugation by $a \in F$:
$$
\tau_a(z)=a z a\inv, \quad z \in F;
$$
it will model the element $a \in F.$
If $\s \in \aut F,$ then
$$
\s \tau_a \s\inv =\tau_{\s(a)}.
$$
Thus, we have interpreted in $F$ the first sort of $\fF.$

Let us prove now that

\begin{Lem}
The set of all primitive elements of $F$
is $\varnothing$-definable in the reduct $\str{\aut
F,F}$ of the structure $\fF.$
\end{Lem}

\begin{pf}
Consider a quasi-conjugation $\f \in \aut F$ with
a canonical basis $\{x\} \cup Y,$ where $C_x=\str{Y}.$
By Proposition \ref{AllWeNeedInCen(Pi)} the set of extremal
involutions $\{\s\}$ of the form
\begin{align*}
\s x &= x\inv, \\
\s y &=x^{2k} y x^{-2k},\quad y \in Y
\end{align*}
is definable in $\aut F$ with the parameter $\f.$
Take one such extremal involution $\s.$
Clearly, $F=\str x * \Fix(\s).$

Let $z$ be a power of a primitive element with $\f z=z\inv.$
Suppose that for each $v\in \Fix(\s),$ the element $zv$ is also
a power of a primitive element. Then we claim that $z=x^{\pm 1}.$
It of course implies that primitive elements are definable.

Indeed, let $z=t^k,$ where $t$ is primitive. Therefore
$(\f t)^k=(t\inv)^k,$ and $\f t=t\inv$
(\cite[I.2.17]{LSch}). By \ref{X_And_X_inv} $t=x$ or
$t=x\inv.$ Thus, $z=x^k$ or $z=x^{-k}$ for some $k \in $ {\bf N}.
Suppose $|k| \ge 2.$ Let $v=u^k,$ where $u$ is a primitive element
in $\Fix(\s).$ Then $zv$ cannot be a power of a primitive element
due to the following result: the equality $a^m b^n c^p=1,$ where $a,b,c$ are
elements of a free group and $m,n,p \ge 2$ implies
that all $a,b,c$ lie in a cyclic subgroup of $F$
(\cite[sect. 6, ch. I]{LSch}).
Really, if $zv=w^p,$ where $w$ is primitive and $|p| >1,$
then both $x$ and $u$ lie in a cyclic subgroup of $F.$
In the case when $p=\pm 1$ we have that $\pm \av w=\av{zv} \in k \Fab$
and $\av w$ is  non-primitive in $\Fab.$ On the other hand, if
$z=x,x\inv,$ then for each $v \in \Fix(\s)$ the element $zv$
is primitive.
\end{pf}

By the Dyer-Scott theorem (see Section \ref{invs}) the
fixed-point subgroup of an involution from $\aut F$ is
a free factor of $F,$ and conversely one can easily
realize an arbitrary free factor of $F$ as the
fixed-point subgroup of some involution. If $\f$ is an
involution in $\aut F$ then
$$
a \in \Fix(\f) \iff \f \tau_a \f\inv =\tau_a.
$$
It also helps us to explain when
involutions $\f,\psi$ have the same fixed-point
subgroups. To model the action of $\aut F$ on the
set of free factors of $F,$ we use the formula
$$
\Fix(\s \f \s\inv)=\s \Fix(\f).
$$

The next is an interpretation of the relation $A=B*C$
on the set of free factors of $F.$ Let us interpret
first the relation $F=A*B.$

(i) Clearly, a free factor $A$ of $F$ has rank one if
and only if $A$ consists of only two distinct
primitive elements. A free factorization
$$
F=A*B
$$
where $\rank A=1$ holds if and only if there is a quasi-conjugation
with a canonical basis $\{x\} \cup Y,$
such that $A=\str x$ and $B=\str Y=C_x^\f.$
Therefore we should write that there is
a quasi-conjugation $\f \in \aut F,$ a primitive element
$z \in F$ with $\f z=z\inv$ such that $z \in A$
and $B=C_z^\f$ (that is $b \in B$ if and only if $\f b=zbz\inv$).

(ii) A soft non-extremal involution $\psi$ with non-trivial
fixed-point subgroup (a member of a definable family
by \ref{SoftsAreDef} and \ref{AllWeNeedInCen(Pi)}),
which satisfies a definable condition
\begin{quote}
`there is a primitive $x \in F$ with $\psi x =x\inv$
such that each primitive $a \in F$ with $\psi
a=a\inv$ is equal to $\psi(w)xw\inv$ for some $w\in F$'
\end{quote}
has by \ref{ConjCriterion} and \ref{TheyMoveToInverse} the
following canonical form
\begin{alignat}2  \label{eqBeads2}
\psi x &=x\inv,   && \\
\psi y &=xyx\inv, &\quad y &\in Y \ne \varnothing \nonumber \\
\psi u &=u,       &u &\in U \ne \varnothing.   \nonumber
\end{alignat}
By \ref{X_And_X_inv} if $a$ is a primitive element,
which $\psi$ sends to its inverse, then $a=v x^{\eps} v\inv,$
where $v \in \Fix(\psi).$ We have
$$
\psi(b)=aba\inv \iff \psi(v\inv b v) =x^\eps (v\inv b v) x^{-\eps}.
$$
Hence $C^\psi_a=v C^\psi_{x^\eps} v\inv.$ Therefore
$$
F =\str a * C_a^\psi * \Fix(\psi).
$$

It follows from the above arguments that the
factorization
$$
F =A*B
$$
where $A$ is of rank at least two, and $B$ is
non-trivial holds if and only if there is an
involution $\psi$ of the form (\theequation) such
that
\begin{itemize}
\item $B$ is equal to the fixed-point subgroup of $\psi;$
\item there is a primitive $a\in F$ with $\psi a=a\inv$
such that $A$ is the least free factor of $F$ which
contains $\str a$ and $C^\psi_a.$
\end{itemize}

Let $A$ be a free proper factor of $F$ of rank $\ge 2.$
Fix a free factor $D$ such that $F=A*D.$ Consider the subgroup
$$
\Sigma=\{\s \in \aut F : \s A=A, \s|_{D} =\id_D\}.
$$
If $\s \in \Sigma$ has the form \eqref{eqBeads2}, then
the restriction $\s$ on $A,$ $\s|_A$ has an analogous canonical
form in $\aut A$ (exactly one block in any canonical basis).
Clearly, $\s|_A$ having such a canonical form
is a quasi-conjugation if and only if $\Fix(\s)=D.$
Therefore $A=B*C,$ where $\rank B=1$ holds if and only if
there are $\s \in \Sigma,$ a primitive $b \in B$ such
that $\s|_A$ is a quasi-conjugation, $b$ and $b\inv$
are only primitive elements in $B,$ and $C=C^\s_b.$
Similarly, one can adapt the method in (ii) above
to construct a definable with the parameter $D$
condition which is equivalent to $A=B*C,$ where
both $B$ and $C$ have rank at least two.
\end{pf}

\section{\it Getting a basis} \label{basis}

Let $F$ be a free group of infinite rank. Consider
free factors $A,B$ of $F$ with bases $\{a_i : i \in I\}$
and $\{b_i : i \in I\},$ respectively such that
$$
F=A*B.
$$
Suppose that $\f$ is an involution with the canonical
form
\begin{alignat*}3
&\f\, &a_i &=a_i\inv, &\quad i &\in I,\\
&\f &b_i &=a_i b_i a_i\inv.
\end{alignat*}
If $\psi$ is an automorphism of $F$ which
fixes all the elements in the subgroup $A$ and commutes
with $\psi,$ then by \ref{TheyMoveToInverse} (ii)
\begin{alignat}2
\psi a_i &=a_i, &\quad i &\in I,\\
\psi b_i &= b_i^{\eps_i}, && \nonumber
\end{alignat}
where $\eps_i=\pm 1.$ Let $\Psi$ denote the set of all automorphisms
of the form (\theequation).

Then $\Psi$ obviously satisfies the following properties:
\begin{itemize}
\item[(a)] each $\psi \in \Psi$ is an involution  and
the fixed-point subgroup of the involution $\psi|_B$
either is trivial, or is a free factor of $B;$

\item[(b)] for every element $b \in B$ there is
$\psi \in \Psi$ such that $\Fix(\psi|_B)$ has finite
rank and $b$ is in $\Fix(\psi|_B);$

\item[(c)] for each  $\psi \in \Psi,$
if $\rank \Fix(\psi|_B) > 1,$ then there are $\psi_1,\psi_2 \in \Psi$ such
that
$$
\Fix(\psi|_B)=\Fix(\psi_1|_B) * \Fix(\psi_2|_B);
$$

\item[(d)] Let $\Psi^1$ denote the set of all elements
in $\Psi$ with $\rank \Fix(\psi|_B)=1;$ then
for every $\psi' \in \Psi^1$ there is a free factor
$C$ of $B$ such that
$$
B=\Fix(\psi'|_B) *C,
$$
and $C$ consists of all the subgroups $\Fix(\psi|_B),$
where $\psi \in \Psi^1 \setminus \{\psi'\}.$
\end{itemize}

The properties (a,b,c,d) imply that the set
$$
\cC=\{b : b \text{ is a primitive element in $\Fix(\psi|_B),$ where } \psi \in 
\Psi^1  \}
$$
is just slightly greater than a `real' basis of $B$:
there exists a basis $\cB$ of $B$ such that $\cC=\cB^\pm.$
Indeed, let $\cB$ is any maximal system of representatives
by an equivalence relation $c \approx d \leftrightarrow (c=d\inv) \vee (c=d)$ on
$\cC.$ Then it follows from (b) and (c)  that $B$ is generated by $\cB,$
because every subgroup of the form $\Fix(\psi|_B)$ of
finite rank can be represented as a free product of fixed-point
subgroups of the rank one of involutions in $\Psi^1|_B.$
The property (d) implies that $\cB$ is a free subset of $B.$

To reduce the set $\cC$ to a basis of $B,$ let us add to
the tuple of parameters $(A,B,\f)$ two involutions
$\pi_0,\pi_1$ such that
\begin{itemize}
\item[(e)] both $\pi_0,\pi_1$ fix $A$ pointwise and preserve $B;$

\item[(f)] $\pi_1|_\cC$ is a bijection between the sets
$$
\cC_0=\{c \in \cC : \pi_0 c=c\inv\} \text{ and } \cC_1=\cC \setminus \cC_0,
$$

\item[(g)] for each $c \in \cC_1$ there is $c_0 \in \cC_0$
such that $\pi_0 c=c_0 c c_0\inv.$
\end{itemize}
Let $\cB_0$ denote the set
$$
\{c \in \cC_0 : (\exists c' \in \cC)\, \pi_0 c'=c c' c\inv\}.
$$
The inverse of each element $c \in \cB_0$
is not in $\cB_0.$ Indeed, if $\pi_0 c'=cc'c\inv,$ then the only primitive
elements $z$ with $\psi_0 z=c\inv z c$ in view of
$$
\pi_0(cc'c\inv)=c\inv(cc'c\inv)c,
$$
are $cc'{}^{\pm1}c\inv,$ but they are not in $\cC.$
Therefore the set
$$
\cB =\cB_0 \cup \pi_1 \cB_0
$$
is a basis of $B.$

It is easy to see that the properties (a,b,c,d,e,f,g)
are definable with parameters $(A,B,\f,\pi_0,\pi_1)$
in the structure $\fF$ over $F$: e.g. $\rank D,$ where
$D$ is an arbitrary free factor of $F,$ is finite if
and only if there is no automorphism in $\aut F$ such
that $\f D$ is a proper subgroup of $D$ (Lemma
\ref{A*B}), we observed above that the condition
$\rank D=1$ is a definable, etc.

To describe by means of first-order logic the tuple of
parameters $(A,B,\f,\pi_0,\pi_1)$ we should explain that
the first two parameters $A$ and $B$ are free factors of
$F$ such that $F=A*B$ and there is an automorphism
$\rho$ of $F$ which maps $A$ onto $B$: $\rho A=B;$ one
can use as the third parameter $\f$ any element
from $\aut F,$ which satisfies the following definable
with the parameters $A,B$ condition: the set of all
automorphisms of $F$
$$
\Psi=\{\psi : \psi\f=\f\psi \text{ and } \psi|_A=\id_A\}
$$
has the properties (a,b,c,d).
There also are no problems in a description of $\pi_0$
and $\pi_1.$

Consider the group
$$
\{\s   : \s B=B \text{ and } \s|_A=\id_A\}
$$
the group $B$ of $F$ as isomorphic
copies of the group $\aut F$ and $F,$ respectively.
Thus, we have proved the following

\begin{Th}
The structure $\str{\aut F,F,\cB}$ {\rm(}with natural relations{\rm)},
where $\cB$ is a basis of $F$ is interpretable in the structure
$\fF$ by means of first-order logic.
\end{Th}

\section{\it Interpretation of set theory} \label{settheo}

\begin{Th} \label{MainTh}
Let $F$ be an infinitely generated free group of rank
$\vk.$ Then the second-order theory of the set $\vk$
and the elementary theory of $\aut F$ are mutually
syntactically interpretable, uniformly in $F$.
\end{Th}

\begin{pf}
It is well-known that the second-order theory
of a set $X$ and first-order theory of the structure
$\str{X,X^X},$ where $X^X$ is the set of all functions
from $X$ to $X,$ are mutually syntactically interpretable.

In the previous section we have interpreted in $\aut F$ the structure
$\str{\aut F,F,\cB},$ where $\cB$ is a basis
of $F.$ Let us partition $\cB$ into two equipotent subsets,
say, $\cB_1$ and $\cB_2,$ by taking two automorphisms
$\pi_0,\pi_1 \in \aut F$ such that
\begin{itemize}
\item $b \in \cB_1 \iff \pi_0 b=b;$
\item $\pi_1 \cB_1=\cB_2.$
\end{itemize}

We shall interpret the set $\vk$ by the set $\cB_1.$
Then we can interpret the set of functions from $\vk$ to
$\vk$ using the set $\Sigma=\{\s\} \subseteq \aut F$ such that
$$
(\forall b \in \cB_1)(\s b=b)\, \& \,
(\forall b \in \cB_1)(\exists b' \in \cB_1)(
 \s(\pi_1 b)=\pi_1 b \cdot b').
$$
Clearly, every $\s \in \Sigma$ determines a $(\pi_0,\pi_1)$-definable
function $b \mapsto b'$ from $\cB_1$ to $\cB_1,$ and, on the other
hand, any function from $\cB_1$ to $\cB_1$ can be coded
in such a way.

Thus, one can interpret (unifomly in $F$) in $\aut F$
the structure $\str{\vk,\vk^\vk}.$

Let $X$ be an infinite set. Let $\strX$ denote
the structure, with the domain $\bigcup_{n \in
\omega} {\cal R}_n(X),$ where $n \in\omega,$ and for
every $n$ ${\cal R}_n(X)$ is the set of all $n$-placed relations on
$X;$ the unique $n$-placed $(n \ge 1)$ basic relation
on $X_{\text{\rm II}}$ says whether
$R(x_1,\ldots,x_{n-1})$ is true or not for any tuple
$x_1,\ldots,x_{n-1} \in X$ and an arbitrary element $R
\in {\cal R}_{n-1}(X).$ Clearly, the elementary
theory of $\strX$ and the second-order theory of $X$
are mutually syntactically interpretable.

It is easy to see that one can interpret in $\strX$
the automorphism group of a free group of rank
$|X|.$ This will complete the proof of the Theorem.

Take a binary operation $f_0 : X \times X \to X$ and
a proper subset $X_0 \subset X$ such that
$|X_0|=|X\setminus X_0|$ (they can be treated as
suitable elements of $\strX$). Then we should explain
that $\str{X;f_0}$ is a group, and any map from $X_0$
to $X$ can be extended to a homomorphism from the
group $\str{X;f_0}$ to a group $\str{X;f},$ where $f$
is an arbitrary group operation on $X.$ Thus, we can interpret
in $\strX$ the structure $\str{F,F^F},$ where $F$ is
a free group of rank $|X|.$ Clearly, the automorphism
group of $F$ is interpretable in the latter structure.
\end{pf}

\begin{Th} \label{ElEquivCrit}
Let $F$ and $F'$ be infinitely generated free
groups of ranks $\vk$ and $\vk',$ respectively.
Then their automorphism groups are elementarily
eqivalent if and only if the cardinals $\vk$ and $\vk'$ are
equivalent in the second-order logic as sets:
$$
\aut F\equiv \aut{F'} \iff \vk \equiv_{\LII} \vk'.
$$
\end{Th}

\begin{pf}
By \ref{MainTh}.
\end{pf}

Theorem \ref{MainTh} implies also that

\begin{Prop}
The first-order theory of the automorphism group
of an infinitely generated free group is undecidable
and unstable.
\end{Prop}


%
\end{document}